\documentclass{article}

\usepackage{amsmath, amssymb, amsthm, amsfonts, bm}
\usepackage{enumerate}
\usepackage{commath}
\usepackage{mathrsfs}
\usepackage[retainorgcmds]{IEEEtrantools}
\usepackage{verbatim}
\numberwithin{equation}{section}
\newcommand{\ud}{\, \mathrm{d}}

\newtheorem{cthm}{Theorem}
\newenvironment{ctheorem}[1]
  {\renewcommand\thecthm{#1}\cthm}
  {\endcthm}

\newtheorem{ccor}{Corollary}

\newtheorem{lemma}{Lemma}[section]
\theoremstyle{definition}
\newtheorem{definition}{Definition}[section]
\theoremstyle{remark}

\newcommand{\Rnum}[1]{\uppercase\expandafter{\romannumeral #1\relax}}
\newcommand{\rnum}[1]{\romannumeral #1\relax}
\newcommand{\mr}[1]{\mathrm{#1}}
\newcommand{\mb}[1]{\mathbb{#1}}
\newcommand{\mc}[1]{\mathcal{#1}}

\DeclareMathOperator{\av}{Av}
\DeclareMathOperator{\mt}{\mathcal{M}_\mathcal{T}}

\title{Sharp integral inequalities for the dyadic maximal operator and applications}
\author{Anastasios D. Delis, Eleftherios N. Nikolidakis}
\date{}

\begin{document}
\maketitle
\footnotetext{ {\em E-mail addresses}: tdelis@math.uoa.gr, lefteris@math.uoc.gr}

\begin{abstract}
We prove a sharp integral inequality for the dyadic maximal function of $\phi\in L^p$. This inequality
connects certain quantities related to integrals of $\phi$ and the dyadic maximal function of $\phi$, under the hypothesis that the variables $\int_X\phi\ud\mu=f,$ $\int_X\phi^q\ud\mu=A,$ $1<q<p,$ are given, where $0<f^q \leq A.$ Additionally, it contains a parameter $\beta>0$ which when it attains a certain value depending only on $f, A, q,$ the inequality becomes sharp. Using this inequality we give an alternative proof of the evaluation of the Bellman 
function related to the dyadic maximal operator of two integral variables.
\end{abstract}

\section{Introduction}\label{sec:1}

It is well known that the dyadic maximal operator on $\mb R^n$ is a useful tool in analysis and is defined by
\begin{equation} \label{eq:1p1}
\mathcal{M}_d \phi(x) = \sup \left\{ \frac{1}{|Q|}\int_Q|\phi(y)|\,\mr dy: x\in Q,\ Q\subseteq \mb R^n\ \text{is a dyadic cube}\right\},
\end{equation}
for every $\phi\in L_\text{loc}^1(\mb R^n)$, where the dyadic cubes are those formed by the grids $2^{-N}\mb Z^n$, for $N=0, 1, 2,\ldots$.
It is also well known that it satisfies the following weak type (1,1) inequality
\begin{equation} \label{eq:1p2}
\left|\left\{x\in\mb R^n : \mathcal{M}_d \phi(x) > \lambda\right\}\right| \leq
\frac{1}{\lambda} \int_{\{\mathcal{M}_d\phi>\lambda\}}|\phi(y)|\,\mr dy,
\end{equation}
for every $\phi\in L^1(\mb R^n)$ and every $\lambda>0$, and which is easily proved to be best possible. Further refinements of \eqref{eq:1p2} can be seen
in \cite{9} and \cite{10}.

Then by using \eqref{eq:1p2} and the well known Doob's method it is not difficult to prove that the following $L^p$ inequality is also true
\begin{equation} \label{eq:1p3}
\|\mathcal{M}_d\phi\|_p \leq \frac{p}{p-1}\|\phi\|_p,
\end{equation}
for every $p>1$ and $\phi\in L^p(\mb R^n)$.
Inequality \eqref{eq:1p3} turns out to be best possible and its sharpness is proved in \cite{16} (for general martingales see \cite{1} and \cite{2}).

One way to study inequalities satisfied by maximal operators is by using the so called Bellman function technique. For example, in order to refine \eqref{eq:1p3} we can insert the $L^1$-norm of $\phi$ as an independent variable in \eqref{eq:1p3}, and try to find the best possible upper bound of $\|\mathcal{M}_d\phi\|_p$, when both the $L^1$ and $L^p$ norms of $\phi$ are given, by evaluating the (Bellman) function of two variables 
\begin{equation} \label{eq:1p4}
B^{(p)}(f,F) = \sup\left\{ \frac{1}{|Q|} \int_Q (\mathcal{M}_d \phi)^p : \phi \geq 0,\ \frac{1}{|Q|}\int_Q\phi=f,\ \frac{1}{|Q|}\int_Q\phi^p=F\right\},
\end{equation}
where $Q$ is a fixed dyadic cube and $f, F$ are such that $0< f^p\leq F$. 

The approach of studying maximal operators by the introduction of the corresponding Bellman function 
was first seen in the work of Nazarov and Treil, \cite{5}, where the authors defined the function  
\begin{multline}\label{eq:1p5}
B_p(f,F,L) =\\ \sup\left\{ \frac{1}{|Q|} \int_Q (\mathcal{M}_d \phi)^p : \frac{1}{|Q|}\int_Q\phi=f,\ \frac{1}{|Q|}\int_Q\phi^p=F, \sup_{R:Q\subseteq R}\frac{1}{|R|}\int_R\phi=L\right\},
\end{multline}
with $p>1$ (as an example they examine the case $p=2$), $Q$ is as above and $\phi$ is non-negative in $L^p(Q),$ $R$ runs over all dyadic cubes containing $Q$ and the variables $F, f, L$ satisfy $0\leq f\leq L,$ $f^p\leq F.$ Exploiting a certain "pseudoconcavity" inequality it satisfies, they construct
the function $4F-4fL+2L^2$ which has the same properties as \eqref{eq:1p5} and provides a good $L^p$ bound for the operator $\mathcal{M}_d$ (see \cite{5} for details).

Both of the above Bellman functions were explicitly computed for the first time by Melas in \cite{3}. 
In fact this was done in the much more general setting of a non-atomic probability space $(X,\mu)$ 
equipped with a tree structure $\mc T$, which is similar to the structure of the dyadic subcubes of 
$[0,1]^n$ (see the definition in Section \ref{sec:2}). Then the associated maximal operator is defined 
as
\begin{equation} \label{eq:1p6}
\mt\phi(x) = \sup\left\{ \frac{1}{\mu(I)}\int_I |\phi|\,\mr d\mu: x\in I\in \mc T\right\},
\end{equation}
for every $\phi\in L^1(X,\mu)$.
Moreover \eqref{eq:1p2} and \eqref{eq:1p3} still hold in this setting and remain sharp. Now if we wish to refine \eqref{eq:1p3} for the general case of a tree $\mathcal{T},$ we should introduce the Bellman function of two variables related to the above maximal operator, which is given by
\begin{equation} \label{eq:1p7}
B_{\mc T}^{(p)}(f,F) = \sup\left\{ \int_X(\mt\phi)^p\,\mr d\mu: \phi\geq 0,\ \int_X\phi\,\mr d\mu = f,\ \int_X\phi^p\,\mr d\mu = F\right\},
\end{equation}
where $0< f^p\leq F$.
This function of course generalizes \eqref{eq:1p4}. In \cite{3} it is proved that
\begin{equation}\label{eq:1p8}
B_{\mc T}^{(p)}(f,F) = F\,\omega_p\!\left(\frac{f^p}{F}\right)^p,
\end{equation}
where $\omega_p : [0,1] \to \bigl[1,\frac{p}{p-1}\bigr]$, is defined by $\omega_p(z) = H_p^{-1}(z)$, 
and $H_p(z)$ is given by $H_p(z) = -(p-1)z^p + pz^{p-1}$. As a consequence $B_{\mc T}^{(p)}(f,F)$ does 
not depend on the tree $\mc T$.
The technique for the evaluation of \eqref{eq:1p7}, that is used in \cite{3}, is based on an effective 
linearization of the dyadic maximal operator that holds on an adequate class of functions called $\mc 
T$-good (see the definition in Section 2), which is enough to describe the problem as is settled in 
\eqref{eq:1p7}. Using this result on suitable subsets of $X$ and several calculus arguments, the 
author also managed to precisely evaluate the corresponding to \eqref{eq:1p5} Bellman function in this context,
\begin{multline}\label{eq:1p9}
B_{\mc T}^{p}(f,F,L) = \sup\bigg\{\int_X (\max(\mathcal{M}_{\mc T} \phi, L)^p \ud \mu :\phi\geq0, \phi \in  L^p(X,\mu),\\ \int_X\phi \ud \mu=f, \int_X\phi^p \ud \mu=F, \bigg\}.
\end{multline}

Now \eqref{eq:1p7} and \eqref{eq:1p9} were computed in \cite{8} in a different way that avoids the calculus arguments involved in \cite{3}. A crucial intermediate result the authors obtain there, in this direction, is the following.
\begin{ctheorem}{A} \label{thm:a}
Let $\phi\in L^p(X,\mu)$ be non-negative, with $\int_X\phi\,\mr d\mu=f$. Then the following inequality is true
\begin{equation} \label{eq:1p10}
\int_X(\mt\phi)^p\,\mr d\mu \leq -\frac{1}{p-1}f^p + \frac{p}{p-1}\int_X \phi\,(\mt\phi)^{p-1}\,\mr d\mu.
\end{equation}
\end{ctheorem} 
The motivation for our work here comes from our wish to refine \eqref{eq:1p7} even further by also considering the $q$-norm, $1<q<p,$ of the function $\phi$ as fixed and to compute the corresponding Bellman function. In particular, our goal was the evaluation of
\begin{multline}\label{eq:1p11}
B_{\mc T}^{p,q}(f,A,F) = \sup\bigg\{\int_X (\mathcal{M}_{\mc T} \phi)^p \ud \mu :\phi\geq0, \phi \in  L^p(X,\mu),\\ \int_X\phi \ud \mu=f, \int_X\phi^q \ud \mu=A \int_X\phi^p \ud \mu=F, \bigg\},
\end{multline}
where $1<q<p,$ and for $f, A, F$ we have $f^q < A < F^{\frac{q}{p}}.$ The new integral variable makes the problem considerably more difficult.

In Sections 3 and 4 we prove our main result, stated in Theorem 1 below. It is an inequality which we 
believe that it can be the corresponding to Theorem \ref{thm:a} intermediate step in the present context, 
towards the evaluation of \eqref{eq:1p11}. Then using this result and entangling a result from \cite{3} 
we prove Corollary \ref{cor:1} below which directly strengthens and generalizes Theorem \ref{thm:a}. 
This will be carried out in Section 5. Finally, also in Section 5, we exploit these results to 
evaluate \eqref{eq:1p7} in a new way.

So our main result is the following.
\begin{ctheorem}{1} \label{thm:1}
Let $q\in(1,p)$, $f>0$ and $\phi\in L^p(X,\mu)$ non-negative, with $\int_X\phi\,\mr d\mu=f.$ Then the inequality
\begin{align}\label{eq:1p12}
\int_X (\mathcal{M}_{\mathcal{T}}\phi)^p \ud \mu & \leq \frac{p(\beta +1)^q}{G(p,q,\beta)}
\int_X (\mathcal{M}_{\mathcal{T}}\phi)^{p-q}\phi^q \ud \mu
+\frac{(p-q)(\beta+1)}{G(p,q,\beta)}f^p \nonumber \\& + \frac{p(q-1)\beta}{G(p,q,\beta)}
f^{p-q}\int_X (\mathcal{M}_{\mathcal{T}}\phi)^q \ud \mu - \frac{p(\beta +1)^q}{G(p,q,\beta)}f^{p-q}\int_X\phi^q\,\mr d\mu,
\end{align}
where $G(p,q,\beta)=p(q-1)\beta+(p-q)(\beta+1),$ is sharp for every $\beta>0.$ If we also assume that
$\int_X\phi^q\,\mr d\mu=A,$ $f^q<A,$ then
\eqref{eq:1p12} is best possible for $\beta=\omega_q(\frac{f^q}{A})-1,$ where $\omega_q$ 
is defined as above, with $q$ in place of $p.$
\end{ctheorem}

Theorem \ref{thm:1}, together with results from \cite{3}, will allow us to prove the following generalization of Theorem \ref{thm:a}, which in turn will lead to the evaluation of \eqref{eq:1p7} in a new way.
\begin{ccor}\label{cor:1}
Let $q\in(1,p)$, $f>0$, and $\phi\in L^p(X,\mu)$ non-negative, with $\int_X\phi\,\mr d\mu=f.$ Then the inequality
\begin{equation}\label{eq:1p13}
\int_X(\mt\phi)^p\,\mr d\mu 
\leq -\frac{q(\beta+1)}{G(p,q,\beta)}f^p +
\frac{p(\beta+1)^q}{G(p,q,\beta)}\int_X (\mathcal{M}_{\mathcal{T}}\phi)^{p-q}\phi^q \ud \mu,
\end{equation}
 is sharp for every $\beta>0,$ where $G(p,q,\beta)$ as above. If we also assume that
$\int_X\phi^q\,\mr d\mu=A,$ $f^q<A,$ then
\eqref{eq:1p13} is best possible for $\beta=\omega_q(\frac{f^q}{A})-1,$
\end{ccor}

We remark here that there are several problems in Harmonic Analysis were Bellman functions arise.
Such problems (including the dyadic Carleson imbedding theorem and weighted inequalities) are described
in \cite{7} (see also \cite{5}, \cite{6}) and also connections to Stochastic Optimal Control are 
provided, from which it follows that the corresponding Bellman functions satisfy certain nonlinear
second-order PDEs. The exact evaluation of a Bellman function is a difficult task which is connected
with the deeper structure of the corresponding Harmonic Analysis problem. Until now several Bellman 
functions have been computed (see \cite{1}, \cite{3}, \cite{5}, \cite{12}, \cite{13}, \cite{14}). The exact computation of \eqref{eq:1p7} has also been given in \cite{11} by L. Slavin,
A. Stokolos and V. Vasyunin, which linked the computation of it to solving certain PDEs of the 
Monge-Amp\`{e}re type, and in this way they obtained an alternative proof of the results in \cite{3} 
for the Bellman function related to the dyadic maximal operator. Also in \cite{15}, using the Monge-
Amp\`{e}re equation approach, a more general Bellman function than the one related to the Carleson 
imbedding theorem has been precisely evaluated thus generalizing the corresponding result in \cite{3}. 
It would be an interesting problem to discover if the Bellman function of three variables defined in 
\eqref{eq:1p11} can be computed using such PDE-based methods.

\section{Preliminaries}\label{sec:2}

In this section we present (without proofs) the background we need from [3], that will be used in all that follows.

Let $(X,\mu)$ be a non-atomic probability space. Two measurable subsets A, B of $X$ will be called
almost disjoint if $\mu(A\cap B)=0$.
\begin{definition}\label{definition:2p1}
A set $\mathcal{T}$ of measurable subsets of $X$ will be called a tree if the following conditions are
satisfied:
\begin{enumerate}
\item[(i)] $X \in \mathcal{T}$ and for every $I \in \mathcal{T}$ we have $\mu(I)>0$.
\item[(ii)]For every $I\in \mathcal{T}$ there corresponds a finite or countable subset
$\mathcal{C}(I) \subseteq \mathcal{T}$ containing at least two elements such that:
\begin{itemize}
\item[(a)]the elements of $\mathcal{C}(I)$ are pairwise almost disjoint subsets of $I$,
\item[(b)]$I=\bigcup \mathcal{C}(I)$.
\end{itemize}
\item[(iii)]$\mathcal{T}= \bigcup_{m\geq 0}\mathcal{T}_{m}$ where $\mathcal{T}_{(0)}=\{X\}$ and
$\mathcal{T}_{(m+1)}=\bigcup_{I\in \mathcal{T}_{(m)}} \mathcal{C}(I)$.
\item[(iv)]We have $\lim_{m \to \infty} \sup_{I \in \mathcal{T}_{(m)}}\mu(I)=0$.
\item[(v)] $\mathcal{T}$ differentiates $L^1(X,\mu)$
\end{enumerate}
This last condition means exactly that the Lesbesgue differentiation theorem holds in the space $L^1(X,\mu)$, with respect to the tree $\mathcal{T}$.

Now we define for any tree $\mathcal{T}$ its exceptional set $E=E(\mathcal{T})$ as follows:
\begin{equation}\label{eq:2p1}
E(\mathcal{T})= \bigcup_{I \in \mathcal{T}} \bigcup_{\substack{J_{1}, J_{2} \in \mathcal{C}(I) \\
J_{1}\neq J_{2}}}
(J_{1}\cap J_{2}).
\end{equation}
It is easy to see that $E(\mathcal{T})$ has measure $0$.
\end{definition}

By induction it can be seen that each family $\mathcal{T}_{(m)}$ consists of pairwise almost disjoint
sets whose union is $X$. Moreover if $x \in X\setminus E(\mathcal{T})$ then for each $m$ there exists
exactly
one $I_{m}(x)$ in $\mathcal{T}_{(m)}$ containing x. For every $m>0$ there is a
$J \in \mathcal{T}_{(m-1)}$ such that $I_{m}(x) \in \mathcal{C}(J)$. Then, since $x \in J,$ we must have
$J=I_{m-1}(x)$, because $x$ does not belong to $E(\mathcal{T})$. Hence the set $\mathcal{A}=\{I \in \mathcal{T}:x \in I\}$ forms a chain
$I_{0}(x)=X\varsupsetneq I_{1}(x)\varsupsetneq \dots$ with $I_{m} \in \mathcal{C}(I_{m-1}(x))$
for every $m>0$. From this remark it follows that if $I, J \in \mathcal{T}$ and
$I \cap J \cap (X\setminus E(\mathcal{T}))$ is nonempty, then $I\subseteq J$ or $J \subseteq I$. In
particular for any $I, J \in \mathcal{T}$, either $\mu(I\cap J)=0$ or one of them is contained in the
other. \\

Given any tree $\mathcal{T}$ we remind that the maximal operator associated to it is defined as follows:
\begin{equation}\label{eq:2p2}
\mathcal{M}_{\mathcal{T}}\phi(x)= \sup \lbrace \frac{1}{\mu(I)} \int_{I} \abs{\phi} \ud \mu \: :\:
x \in I \in \mathcal{T} \rbrace,
\end{equation}
for every $\phi \in L^1(X,\mu)$.

Next we describe the linearization procedure for the operator $\mathcal{M}_{\mathcal{T}}$.
Let $\phi \in L^1(X,\mu)$ be a nonnegative function and for any $I \in \mathcal{T}$ let
\begin{equation}\label{eq:2p3}
\av_I (\phi) = \frac{1}{\mu(I)}\int_{I} \phi \ud \mu.
\end{equation}
We will say that $\phi$ is $\mathcal{T}$-good if the set
\begin{equation}\label{eq:2p4}
\Lambda_{\phi}=\lbrace x \in X\setminus E(\mathcal{T}) \: : \: \mathcal{M}_{\mathcal{T}}\phi(x)
>\av_I(\phi) \; for \; all \; I \; \in \mathcal{T} \; such \; that \; x \in I \rbrace
\end{equation}
has $\mu$-measure zero.

For any such function and every $x \in X\setminus (E(\mathcal{T})\cup \Lambda_{\phi})$ (i.e. for
$\mu$-almost every $x$ in X) we define $I_{\phi}(x)$ to be the largest element in the nonempty set
$\{ I \in \mathcal{T}: x \in I \; \text{and} \; \mathcal{M}_{\mathcal{T}}\phi(x)=\av_I(\phi)\}$.

Also given any $I \in \mathcal{T}$ let
\begin{equation}\label{eq:2p5}
A(\phi, I) = \{x \in X\setminus (E(\mathcal{T})\cup \Lambda_{\phi}) : I_{\phi}(x)=I\}
\subseteq I
\end{equation}
and
\begin{equation}\label{eq:2p6}
\mathcal{S}_{\phi}= \{I \in \mathcal{T}: \mu(A(\phi, I)>0\} \cup \{X\}.
\end{equation}
It is clear that
\begin{equation}\label{eq:2p7}
\mathcal{M}_{\mathcal{T}}\phi = \sum_{I \in \mathcal{S}_{\phi}}Av_I(\phi)\chi_{A(\phi,I)},
\; \text{almost everywhere},
\end{equation}
where $\chi_B$ denotes the characteristic function of $B\subset X.$
Now we define the correspondence $I \to I^*$ with respect to $\mathcal{S}_{\phi}$ for $I \neq X$
in the following manner:
$I^*$ is the smallest element of $\{J \in \mathcal{S}_{\phi}: \:I \subsetneq J \}$.

It is clear that the sets $A_I=A(\phi, I), \; I \in \mathcal{S}_{\phi},$ are pairwise disjoint and
since $\mu(\cup_{J \notin \mathcal{S}_{\phi}}A_J )=0$ their union has full measure.

In the following Lemma we present several important properties of the sets defined above. At this point we define two
measurable sets $A$ and $B$ to be almost equal if $\mu(A \setminus B)=\mu(B\setminus A)=0$ and in this    
case we write $A \thickapprox B$

\begin{lemma}\label{lem:2p2}
\begin{enumerate}
\item[(i)]If $I,J \in \mathcal{S}_{\phi}$ then either $A_J\cap I = \varnothing$ or
$J \subseteq I$.
\item[(ii)] If $I \in \mathcal{S}_{\phi}$ then there exists $J \in \mathcal{C}(I)$ such that
$J \notin \mathcal{S}_{\phi}$
\item[(iii)]For every $I \in \mathcal{S}_{\phi}$ we have $I \thickapprox \bigcup_
{\mathcal{S}_{\phi} \ni J \subseteq I} A_J$.
\item[(iv)]For every $I \in \mathcal{S}_{\phi}$ we have $A_I \thickapprox
I \setminus \bigcup_{J \in \mathcal{S}_{\phi}: \; J^*=I}J$ and so
\begin{equation}\label{eq:2p8}
\mu(A_I)= \mu(I) - \sum_{J \in \mathcal{S}_{\phi}: \; J^*=I}\mu(J).
\end{equation}
\end{enumerate}
\end{lemma}

From the above we get
\begin{equation}\label{eq:2p9}
\av_I(\phi) = \frac{1}{\mu(I)}\sum_{J \in \mathcal{S}_{\phi} : \; J \subseteq I}
\int_{A_J}\phi \ud \mu .
\end{equation}
Now we fix $q>1$. Following \cite{3} we set
\begin{equation}\label{eq:2p10}
x_I = a_I^{-1+\frac{1}{q}} \int_{A_I}\phi \ud \mu
\end{equation}
for every $I \in \mathcal{S}_{\phi}$ where $a_I = \mu(A_I)$ (in case where
$\mu(A_X)=0$ we set $x_X=0$) and from H$\ddot{o}$lder's inequality and Lemma \ref{lem:2p2} we get
\begin{equation}\label{eq:2p11}
\mathcal{M}_{\mathcal{T}}\phi= \sum_{I \in \mathcal{S}_{\phi}}\left(\frac{1}{\mu(I)}
\sum_{J \in \mathcal{S}_{\phi} : \; J \subseteq I} a_J^{1/\acute{q}}x_J
\right) \chi_{A_I}
\end{equation}
$\mu$-almost everywhere, where $\acute{q}=\frac{q}{q-1}$ is the dual exponent of q, and also
\begin{equation}\label{eq:2p12}
\int_X \phi^q \ud \mu= \sum_{I \in \mathcal{S}_{\phi}} \int_{A_I} \phi^q \ud \mu
\geq \sum_{I \in \mathcal{S}_{\phi}}x_I^q.
\end{equation}
So we have
\begin{equation}\label{eq:2p13}
\int_X (\mathcal{M}_{\mathcal{T}}\phi)^q \ud\mu= \sum_{I \in \mathcal{S}_{\phi}}
\left(\frac{1}{\mu(I)} \sum_{J \in \mathcal{S}_{\phi} : \; J \subseteq I}
a_J^{1/\acute{q}}x_J \right)^q a_I = \sum_{I \in \mathcal{S}_{\phi}}a_I y_I^q
\end{equation}
where 
\begin{equation}\label{eq:2p14}
y_I=\av_I(\phi)=\frac{1}{\mu(I)}\sum_{J \in \mathcal{S}_{\phi} : \; J \subseteq I}a_J^{1/\acute{q}}x_J.
\end{equation}

\section{Proof of \eqref{eq:1p12}}\label{sec:3}

We shall first prove \eqref{eq:1p12} for the class of $\mathcal{T}$-good functions.
Let $\phi: (X,\mu)\to \mb R^+$ be $\mc T$-good and such that $\int_X\phi\,\mr d\mu=f$ and $\int_X\phi^q\,\mr d\mu=A$.
We use the linearization technique mentioned in Section 2. From \eqref{eq:2p11} and \eqref{eq:2p14}, if we set
\begin{equation*}
F'=\int_X(\mt \phi)^{p-q}\phi^q\ud\mu,
\end{equation*}
we get
\begin{align} \label{eq:3p1}
F' &
=\int_X \sum_{I\in \mathcal{S}}y_I^{p-q} \chi_{A_I}\phi^q \ud \mu
=\sum_{I\in \mathcal{S}}y_I^{p-q}\int_{A_I} \phi^q \ud \mu \nonumber \\&
=\sum_{\substack{I\in \mathcal{S}\\I \neq X}}y_I^{p-q}\int_{A_I} \phi^q \ud \mu +y_X^{p-q}A -
y_X^{p-q}\sum_{\substack{I\in \mathcal{S}\\I^*= X}}\int_I \phi^q \ud \mu,
\end{align}
where for the last equality in \eqref{eq:3p1} we used Lemma \ref{lem:2p1}(\rnum 4).
Lemma \ref{lem:2p1}(\rnum 3) and the definition of the correspondence $I\to I^*$ imply
\begin{equation}\label{eq:3p2}
\sum_{\substack{I\in \mathcal{S}\\I^*= X}}\int_I \phi^q \ud \mu
=\sum_{\substack{I\in \mathcal{S}\\I^*= X}}\sum_{\substack{J\in \mathcal{S}\\J \subseteq I}}
\int_{A_J} \phi^q \ud \mu
=\sum_{\substack{I\in \mathcal{S}\\I\neq X}}\int_{A_I} \phi^q \ud \mu.
\end{equation}
Moreover, it is easy to see that
$$x_I^q= a_I^{q-1}(y_I\mu(I) -\sum_{\substack{J\in \mathcal{S}\\J^*=I}}y_J\mu(J))^q$$ and
$$\int_{A_I}\phi^q\ud\mu \geq x_I^q.$$
So using H\"{o}lder's inequality in the form
\begin{equation} \label{eq:3p3}
\frac{(\lambda_1 + \lambda_2 + \ldots + \lambda_m)^q}{(\sigma_1 + \sigma_2 + \ldots + \sigma_m)^{q-1}} \leq \frac{\lambda_1^q}{\sigma_1^{q-1}} + \frac{\lambda_2^q}{\sigma_2^{q-1}} + \ldots + \frac{\lambda_m^q}{\sigma_m^{q-1}},
\end{equation}
which holds for every $\lambda_i\geq 0,\ \sigma_i>0$, since $q>1,$ Lemma \ref{lem:2p1}(\rnum 4) and the properties of the correspondence $I\to I^*,$ \eqref{eq:3p1} becomes
\begin{align}\label{eq:3p4}
F'&=\sum_{\substack{I\in \mathcal{S}\\I \neq X}}y_I^{p-q}\int_{A_I} \phi^q \ud \mu +y_X^{p-q}A- y_X^{p-q}\sum_{\substack{I\in \mathcal{S}\\I\neq X}}\int_{A_I} \phi^q \ud \mu \nonumber
\\
&=\sum_{I\in \mathcal{S}}(y_I^{p-q} - y_X^{p-q})\int_{A_I} \phi^q \ud \mu +y_X^{p-q}A \nonumber\\
&\geq \sum_{I\in \mathcal{S}}(y_I^{p-q} - y_X^{p-q})x_I^q +y_X^{p-q}A \nonumber \\
&= \sum_{I\in \mathcal{S}}(y_I^{p-q} - y_X^{p-q})
\frac{(y_I\mu(I) -\sum_{\substack{J\in \mathcal{S}\\J^*=I}}y_J\mu(J))^q}
{(\mu(I)-\sum_{\substack{J\in \mathcal{S}\\J^*=I}}\mu(J))^{q-1}} +y_X^{p-q}A\nonumber \\
&\geq \sum_{I\in \mathcal{S}}(y_I^{p-q} - y_X^{p-q})
\Big(\frac{(y_I\mu(I))^q}{(\tau_I\mu(I))^{q-1}}-\sum_{\substack{J\in\mathcal{S}\\J^*=I}}
\frac{(y_J\mu(J))^q}{((\beta+1)\mu(J))^{q-1}}\Big) +y_X^{p-q}A\nonumber \\
&= K -\sum_{\substack{I\in \mathcal{S}\\I \neq X}}(y_{I^*}^{p-q} - y_X^{p-q})
\frac{(y_I\mu(I))^q}{((\beta+1)\mu(I))^{q-1}}\notag \\&=
K -\sum_{\substack{I\in \mathcal{S}\\I \neq X}}y_{I^*}^{p-q}y_I^q\frac{\mu(I)^q}{((\beta+1)\mu(I))^{q-1}} + y_X^{p-q}\sum_{\substack{I\in \mathcal{S}\\I \neq X}}\frac{(y_I\mu(I))^q}{((\beta+1)\mu(I))^{q-1}}
\end{align}
provided that the $\tau_I>0$ satisfy
$\tau_I\mu(I)-(\beta+1)\sum_{J^*=I}\mu(J))=\mu(I)-\sum_{J^*=I}\mu(J)$, which in turn gives
\begin{equation}\label{eq:3p5}
\tau_I=\beta+1-\beta \rho_I,
\end{equation}
with $\rho_I=\frac{a_I}{\mu(I)},$ and
\begin{equation}\label{eq:3p6}
K=\sum_{I\in \mathcal{S}}(y_I^{p-q} - y_X^{p-q})\frac{(y_I\mu(I))^q}{(\tau_I\mu(I))^{q-1}}
+y_X^{p-q}A,
\end{equation}
We now use the following elementary inequality,
\[
p x^q\!\cdot\! y^{p-q} \leq q x^p + (p-q) y^p,
\]
which holds since $1< q < p$, for any $x, y>0,$ to get
\begin{equation}\label{eq:3p7}
F' \geq K-\frac{p-q}{p}\sum_{\substack{I\in \mathcal{S}\\I \neq X}}y_{I^*}^p\frac{\mu(I)}{(\beta+1)^{q-1}}
-\frac{q}{p}\sum_{\substack{I\in \mathcal{S}\\I \neq X}}y_{I}^p\frac{\mu(I)}{(\beta+1)^{q-1}}
+y_X^{p-q}\sum_{\substack{I\in \mathcal{S}\\I \neq X}}\frac{y_I^q\mu(I)}{(\beta+1)^{q-1}}.
\end{equation}
From Lemma \ref{lem:2p2} (\rnum 4),
\begin{align} \label{eq:3p8}
\sum_{\substack{I\in S\\ I\neq X}} \mu(I)y_{I^*}^p &=\sum_{I\in \mathcal{S}}\sum_{\substack{J\in S \\J^*= I}}\mu(J)y_I^p= \sum_{I\in S}(\mu(I)-a_I)y_I^p \nonumber \\&
=y_X^p +\sum_{\substack{I\in S\\ I\neq X}}\mu(I)y_I^p - \sum_{I\in S}a_I y_I^p.
\end{align}
So 
\begin{align}\label{eq:3p9}
F' &\geq
K-\sum_{\substack{I\in \mathcal{S}\\I \neq X}}(y_I^{p-q} - y_X^{p-q})
\frac{y_{I}^q\mu(I)}{(\beta+1)^{q-1}}-\frac{(p-q)y_X^p}{p(\beta+1)^{q-1}}
+\frac{p-q}{p}\sum_{I\in \mathcal{S}}\frac{a_Iy_{I}^p}{(\beta+1)^{q-1}} \nonumber \\
&= \sum_{I\in \mathcal{S}}(y_I^{p-q} - y_X^{p-q})\frac{1}{\rho_I}
\left(\frac{1}{(\beta+1-\beta\rho_I)^{q-1}}-\frac{1}{(\beta+1)^{q-1}}\right)a_Iy_{I}^q\nonumber \\
&-\frac{p-q}{p}\frac{y_X^p}{(\beta+1)^{q-1}}
+\frac{p-q}{p}\sum_{I\in \mathcal{S}}\frac{a_Iy_{I}^p}{(\beta+1)^{q-1}} +y_X^{p-q}A,
\end{align}
after we have expanded $K.$
Note now that
\[
\frac{1}{(\beta+1-\beta x)^{q-1}} - \frac{1}{(\beta+1)^{q-1}} \geq \frac{(q-1)\beta x}{(\beta+1)^q},
\]
by the mean value theorem on derivatives for all $x\in [0,1]$, so \eqref{eq:3p9} becomes
\begin{align} \label{eq:3p10}
F' &\geq \left(\frac{(q-1)\beta}{(\beta+1)^{q}}+ \frac{p-q}{p(\beta+1)^{q-1}}\right)
\int_X (\mathcal{M}_{\mathcal{T}}\phi)^p \ud \mu \nonumber \\&
-\frac{p-q}{p}\frac{f^p}{(\beta+1)^{q-1}}
-f^{p-q}\frac{(q-1)\beta}{(\beta+1)^{q}}\int_X(\mathcal{M}_{\mathcal{T}}\phi)^q \ud \mu
+f^{p-q}A
\end{align}
for every $\beta > 0$. Rearranging the terms, we get \eqref{eq:1p12} for $\mc{T}$-
good functions.

For the general $\phi \in L^p(X,\mu)$ with $\int_X \phi=f$ and $\int_X \phi^q=A$, $f^q<A$, $1<q<p,$ \eqref{eq:1p12} is proved as follows. We consider the sequence
$\{\phi_m\}$, where $\phi_m = \sum_{I\in \mathcal{T}_{(m)}}\av_I(\phi)\chi_I$ and we set
\begin{equation*}
\Phi_m= \sum_{I\in \mathcal{T}_{(m)}}\max\{\av_I(\phi):\; I\subseteq J \in \mathcal{T}\}\chi_I=
\mathcal{M}_\mathcal{T}\phi_m,
\end{equation*}
since $\av_J(\phi)=\av_J(\phi_m)$ whenever $I\subseteq J \in \mathcal{T}$. It is easy to see that
\begin{equation}\label{eq:3p11}
\int_X \phi_m \ud\mu=\int_X \phi \ud\mu=f, \qquad
\int_X \phi_m^q \ud\mu \leq \int_X \phi^q \ud\mu
\end{equation}
for all m and that $\Phi_m$ converges monotonically almost everywhere
to $\mathcal{M}_\mathcal{T}\phi$. Since $\phi_m$ is easily seen to be $\mc{T}$-good, from what we have 
just shown,
\begin{align}\label{eq:3p12}
\int_X \Phi_m^p \ud \mu & \leq \frac{p(\beta +1)^q}{G(p,q,\beta)}
\int_X \Phi_m^{p-q}\phi_m^q \ud \mu +\frac{(p-q)(\beta+1)}{G(p,q,\beta)}f^p \nonumber \\& + \frac{p(q-1)\beta}{G(p,q,\beta)}
f^{p-q}\int_X \Phi_m^q \ud \mu - \frac{p(\beta +1)^q}{G(p,q,\beta)}f^{p-q}\int_X \phi_m^q.
\end{align}
Since $\mathcal{T}$ differentiates $L^1(X,\mu)$
and by the definition of $ \phi_m,$ if
$\{I_m(x)\}$ is the chain of elements of $\mathcal{T}$ which contain $x\in X,$ then
\begin{equation}\label{eq:3p13}
\lim_{m \to \infty}\phi_m(x) = \lim_{m \to \infty}\av_{I_m(x)}(\phi) = \phi(x)
\end{equation}
and $\phi_m \leq \Phi_m$
Taking limits using the monotone and dominated convergence theorems and Fatou's lemma, we
obtain \eqref{eq:1p12} for the general $\phi \in L^p(X,\mu).$

\section{Proof of Theorem \ref{thm:1}}\label{sec:4}

We now move on to show that \eqref{eq:1p12} is sharp. To do this we shall use a result from \cite{8}
stated in Theorem \ref{thm:2} below. What makes it particularly useful in our case is that it is valid 
for any functions $G_1, G_2$ satisfying the specific properties mentioned. We remind here that the 
decreasing rearrangement $\phi^*:(0, \infty) \to [0, \infty]$ of a measurable function $\phi :X \to 
\mathbb{R},$ is defined as 
$$\phi^*(t)=\inf\{s: d_\phi(s) \leq t\},$$ with $d_\phi$ the distribution function of $\phi.$
\begin{ctheorem}{2}\label{thm:2}
The following is true
\begin{align}\label{eq:4p1}
&\sup\big\{\int_KG_1(\mt\phi)G_2(\phi)\ud\mu : \phi^*=g,\text{ } K \subseteq X \text{ measurable, with } \mu(K)=k\big\}\nonumber \\& = \int_0^kG_1(\frac{1}{t}\int_0^tg)G_2(g(t))\ud t,
\end{align}
where $G_i:[0,+\infty] \to [0,+\infty],$ $i=1, 2,$ are increasing functions, $g:(0, 1] \to 
\mathbb{R}$ is non-increasing and $\phi^*$ is the decreasing rearrangement of the function $\phi.$
\end{ctheorem}
So, with $X$ in place of $K$ and from well known properties of the decreasing rearrangement, it is now easy to see that 
\begin{align}\label{eq:4p2}
&\sup \big\{\int_XG_1(\mt\phi)G_2(\phi)\ud\mu : \phi \geq 0, \text{ measurable, with } 
\int_X\phi=f\big\}\nonumber \\& =\sup\big\{\int_0^1G_1(\frac{1}
{t}\int_0^tg)G_2(g(t))\ud t  : g:(0, 1] \to \mathbb{R}, \text{ non-increasing, } \int_0^1g=f\big\}.
\end{align}
 
Let $\beta>0$ and define $g:(0, 1] \to \mathbb{R}$ by
\begin{equation}\label{eq:4p3}
g(t)=\frac{f}{\beta+1}t^{-1+\frac{1}{\beta+1}}.
\end{equation}
It is easy to see that $\int_0^1g=f,$ for every $\beta>0,$ $\frac{1}{t}\int_0^tg=(\beta+1)g(t)$ for every $t\in(0,1]$ and, after straightforward calculations,
that
\begin{align}\label{eq:4p4}
\int_0^1 (\frac{1}{t}&\int_0^tg)^p \ud \mu  = \frac{p(\beta +1)^q}{G(p,q,\beta)}
\int_0^1(\frac{1}{t}\int_0^tg)^{p-q}g^q \ud \mu
+\frac{(p-q)(\beta+1)}{G(p,q,\beta)}f^p \nonumber \\& + \frac{p(q-1)\beta}{G(p,q,\beta)}
f^{p-q}\int_0^1 (\frac{1}{t}\int_0^tg)^q \ud \mu - \frac{p(\beta +1)^q}{G(p,q,\beta)}f^{p-q}\int_0^1g^q\,\mr d\mu.
\end{align}
This, together with \eqref{eq:4p2}, proves the sharpness of \eqref{eq:1p12} for the 
first case stated in Theorem \ref{thm:1}. Now, \eqref{eq:4p2} is valid if we add $\int_X\phi^q=A$ in the brackets of the left hand side of \eqref{eq:4p2} and $\int_0^1g^q=A$ to the right, that is if we consider the q-norms of
the corresponding functions as fixed. So in case $\int_X\phi^q \ud\mu=A,$ all we need to do is choose 
$\beta>0$ so, that $\int_0^1g^q=A,$ with $g$ as in \eqref{eq:4p3}. The appropriate value is easily seen to be the one given in Theorem \ref{thm:1} and the proof is complete.

\section{Applications}\label{sec:5}

\begin{proof}[Proof of Corollary \ref{cor:1}]~\\
Since $\int_X \phi \ud \mu =f$, from (4.25) in [3], we know that
\begin{equation}\label{eq:5p1}
\int_X(\mathcal{M}_{\mathcal{T}}\phi)^q \ud \mu \leq \frac{\beta+1}{\beta}
\frac{(\beta+1)^{q-1}\int_X \phi^q \ud \mu-f^q}{q-1}
\end{equation}
for every $\beta>0,$ for $\phi$ a $\mc{T}$-good function.
Plugging this into \eqref{eq:3p10} we get \eqref{eq:1p13} for $\mc{T}$-good functions and defining
$\phi_m$ and $\Phi_m$ as in Section 3, we get \eqref{eq:1p13} for the general $\phi \in L^p(X,\mu),$ using the monotone convergence theorem. Sharpness is proved for both cases in the same way it has been
proved for \eqref{eq:1p12}. We only need to observe that with $g$ as in \eqref{eq:4p3}
\begin{equation}\label{eq:5p2}
\int_0^1 (\frac{1}{t}\int_0^tg)^p \ud \mu =  -\frac{q(\beta+1)}{G(p,q,\beta)}f^p +
\frac{p(\beta+1)^q}{G(p,q,\beta)}\int_0^1(\frac{1}{t}\int_0^tg)^{p-q}g^q \ud \mu.
\end{equation}
\end{proof}

Our final application is to derive the least upper bound for $\int_X\left(\mt\phi\right)^p\mr d\mu$, when on $\phi$ we impose the conditions
$\int_X\phi\mr d\mu=f$ and $\int_X\phi^p\mr d\mu=F$ (where $f,F$ are fixed, satisfying $0\leq f^p\leq F$), by using the proof of inequality \eqref{eq:1p13},
for an arbitrary $q$ belonging to $(1,p)$, and a suitable value of $\beta$, depending on $q, p, f$ and $F$. That is we find the main Bellman function of two variables, \eqref{eq:1p7}, associated to the dyadic maximal operator. We proceed to this as follows.

Fix $q\in (1,p)$. First of all it is easy to see that for the above $f, F$, there exists $\beta \in (0,\frac{1}{p-1})$, such that

\begin{equation}\label{eq:5p3}
h_{\beta}(\beta+1)F=\frac{q}{p}\frac{1}{(\beta+1)^{q-1}}f^p,
\end{equation}
where $h_{\beta}(y)$ is defined, for every $y>1$, by $h_{\beta}(y)=y^{p-q}-A_{\beta}y^p$ and
$A_{\beta}$ is defined by
 
\begin{equation}\label{eq:5p4}
A_{\beta}=\frac{(q-1)\beta}{(\beta+1)^q}+\frac{p-q}{p}\frac{1}{(\beta+1)^{q-1}}.
\end{equation}
For this existence, we just need to define the function
$$G(\beta)=\frac{1}{(\beta+1)^{p-1}[1-\beta(p-1)]},$$
of $\beta\in (0,\frac{1}{p-1})$, and note that $G(0+)=1$ and $G(\frac{1}{p-1}-)=+\infty$. Thus
there exists $\beta\in (0,\frac{1}{p-1})$, such that $G(\beta)=\frac{F}{f^p}\geq1$. If this last condition is true we easily see, after some simple calculations,
that $h_{\beta}(\beta+1)F=\frac{q}{p}\frac{1}{(\beta+1)^{q-1}}f^p$, which is \eqref{eq:5p3}.

Now, because of \eqref{eq:1p13}, for any $\phi \in L^p(X,\mu),$ and for this value of $\beta$, the following inequality holds
\[
\int_X\phi^q\left(\mt\phi\right)^{p-q}\mr d\mu\geq A_{\beta}\int_X\left(\mt\phi\right)^p\mr d\mu + \frac{q}{p}\frac{1}{(\beta+1)^{q-1}}f^p. 
\]
Applying H\"{o}lder's inequality on the left side of the above inequality we obtain
\[
F^{q/p}\big(\int_X\left(\mt\phi\right)^p\mr d\mu\big)^{(p-q)/p}\geq A_{\beta}\int_X\left(\mt\phi\right)^p\mr d\mu+\frac{q}{p}\frac{1}{(\beta+1)^{q-1}}f^p
\]
or equivalently, by dividing both sides by $F$,
\[
I_{\phi}^{(p-q)/p}\geq A_{\beta}I_{\phi}+\frac{q}{p}\frac{1}{(\beta+1)^{q-1}}\frac{f^p}{F},
\]
where in the last inequality we denote $I_{\phi}=\frac{\int_X\left(\mt\phi\right)^p\mr d\mu}{F}$, which in turn means that
\begin{equation}\label{eq:5p5}
h_{\beta}(I_{\phi}^{1/p})\geq \frac{q}{p}\frac{1}{(\beta+1)^{q-1}}\frac{f^p}{F}.
\end{equation}

Now for any $\beta\in (0,\frac{1}{p-1})$, we prove that the function $h_{\beta}$, with domain $(1,+\infty)$, is strictly decreasing. For this proof
we proceed in the following way. We have that $\frac{d}{dy}h_{\beta}(y)=y^{p-1}[(p-q)y^{-q}-pA_{\beta}]<y^{p-1}[(p-q)-pA_{\beta}]$, where the inequality in the
last relation is true due to the fact that $y$ is greater than $1$.
Now we claim that $A_{\beta}>\frac{p-q}{p}$, for any $q\in[1,p]$ and $\beta\in (0,\frac{1}{p-1})$. For this reason, we consider $A_{\beta}$ as a function
of $\beta$, in the above mentioned domain and denote it as $K(\beta)$. Then $K(0)=\frac{p-q}{p}$, so we just need to prove that $K(\beta)$ is strictly
increasing. For this purpose we evaluate $\frac{d}{d\beta}K(\beta)$, which as can be easily seen by using \eqref{eq:5p4} is equal to $\frac{(q-1)q[1-\beta(p-1)]}{p(\beta+1)^{q+1}}$,
which is positive for any $\beta$ as above.
By the above discussion we conclude that $\frac{d}{dy}h_{\beta}(y)<0$, for any $y>1$. 

Thus from \eqref{eq:5p5} we have as a consequence that
$I_{\phi}^{1/p}\leq h_{\beta}^{-1}(L)$, where $L=\frac{q}{p}\frac{1}{(\beta+1)^{q-1}}\frac{f^p}{F}$. This conclusion holds, if we suppose that
$I_{\phi}>1$, which may be assumed, since in the opposite case we have nothing to prove. We finally reach the inequality
 
\begin{equation}\label{eq:5p6}
\int_X\left(\mt\phi\right)^p\mr d\mu\leq F(h_{\beta}^{-1}(L))^p
\end{equation}

Having now in mind that \eqref{eq:5p3} holds, we show that $h_{\beta}^{-1}(L)=\omega_p\!\left(\frac{f^p}{F}\right)$,
where $\omega_p$ is defined in the Introduction.
Indeed, by \eqref{eq:5p3}, we immediately conclude that $h_{\beta}^{-1}(L)=\beta+1$, so we just need to prove that
$\beta+1=\omega_p\!\left(\frac{f^p}{F}\right)$. Equivalently this means that $H_p(\beta+1)=\frac{f^p}{F}$. But by \eqref{eq:5p3},
we easily see that
$$\frac{p}{q}(\beta+1)^{q-1}[(\beta+1)^{p-q}-A_{\beta}(\beta+1)^p]=\frac{f^p}{F}.$$
After simple calculations in the left side of the above equality, the real number $q$ is cancelled giving us the quantity
$$-(p-1)(\beta+1)^p+ p(\beta+1)^{p-1},$$
which is exactly $H_p(\beta+1)$. In this way we derive that
 
\[
B_{\mc T}^{(p)}(f,F) \leq F\,\omega_p\!\left(\frac{f^p}{F}\right)^p.
\]
This establishes the least upper bound we need to find for the quantity of interest for the general $\phi\in L^p(X,\mu)$.
Note finally that the opposite inequality is also true, as can be concluded immediately by the sharpness of inequality \eqref{eq:1p13},
which is best possible for any fixed values of $f$ and $\beta.$
Thus we have equality in the above inequality, and our evaluation of the Bellman function of two variables for the dyadic maximal operator is completed.

Anastasios D. Delis, Eleftherios N. Nikolidakis, National and Kapodistrian University of Athens, Department of Mathematics, Panepistimioupolis, Zografou 157 84,
Athens, Greece.


\begin{thebibliography}{99}
	






\bibitem{1}
	D. L. Burkholder,
	\emph{Martingales and Fourier Analysis in Banach spaces},
	C.I.M.E. Lectures, Varenna, Como, Italy, 1985,
	Lecture Notes Math. 1206 (1986),
	81--108.
	
\bibitem{2}
	D. L. Burkholder,
	\emph{Explorations in martingale theory and its applications},
	\'{E}cole d' \'{E}t\'{e} de Probabiliti\'{e}s de Saint-Flour XIX--1989,
	Lecture Notes Math. 1464 (1991),
	1--66.
	
\bibitem{3}
	A. D. Melas,
	\emph{The Bellman functions of dyadic-like maximal operators and related inequalities},
	Adv. in Math. 192 (2005),
	310--340.

\bibitem{4}
	A. D. Melas,
	\emph{Sharp general local estimates for dyadic-like maximal operators and related Bellman functions},
	Adv. in Math. 220 (2009),No 2
	367--426.
	
\bibitem{5}
	F. Nazarov, S. Treil,
	\emph{The hunt for a Bellman function: Applications to estimates for singular integral operators and to other classical problems of harmonic analysis},
	St. Petersburg Math. J. 8 no. 5 (1997),
	721--824

\bibitem{6}
	F. Nazarov, S. Treil and A. Volberg,
	\emph{The Bellman functions and two-weight inequalities for Haar multipliers},
	Journ. Amer. Math. Soc. 12 no. 4 (1999),
	909--928.
	
\bibitem{7}
	F. Nazarov, S. Treil and A. Volberg,
	\emph{Bellman function in stochastic optimal control and harmonic analysis (how our Bellman function got its name)},
	Oper. Theory: Advances and Appl. 129 (2001), 393-424, Birkhauser, Verlag. MR1882704 (2003b:49024).	

\bibitem{8}
	E. N. Nikolidakis, A. D. Melas,
	\emph{A sharp integral rearrangement inequality for the dyadic maximal operator and applications},
	Appl. and  Comp. Harmonic Anal., 38 (2015), Issue 2,
    242--261.

\bibitem{9}
	E. Nikolidakis,
	\emph{Optimal weak type estimates for dyadic-like maximal operators},
	Ann. Acad. Scient. Fenn. Math. 38 (2013),
	229--244.
	
\bibitem{10}
	E. Nikolidakis,
	\emph{Sharp weak type inequalities for the dyadic maximal operator},
	J. Fourier. Anal. Appl., 19 (2012),
	115--139.

\bibitem{11}
	L. Slavin, A. Stokolos, V. Vasyunin,
	\emph{Monge-Amp\`{e}re equations and Bellman functions: The dyadic maximal operator}
	C. R. Math. Acad. Sci. Paris S\'{e}r. I. 346 (2008),
	585--588.
	
\bibitem{12}
	L. Slavin, A. Volberg,
	\emph{The explicit BF for a dyadic Chang-Wilson-Wolff theorem. The $s$-function and the exponential integral},
	Contemp. Math. 444.
	Amer. Math. Soc., Providence, RI, 2007.
	
\bibitem{13}
	V. Vasyunin,
	\emph{The sharp constant in the reverse H\"{o}lder inequality for Muckenhoupt weights},
	St. Petersburg Math. J., 15 (2004), no. 1,
	49--75.
	
\bibitem{14}
	V. Vasyunin, A. Volberg,
	\emph{The Bellman functions for a certain two weight inequality: The case study},
	St. Petersburg Math. J., 18 (2007), No. 2,
	p 201--222.
	
\bibitem{15}
	V. Vasyunin, A. Volberg,
	\emph{Monge-Amp\`{e}re equation and Bellman optimization of Carleson embedding theorems, Linear and complex analysis},
	195--238, Amer. Math. Soc. Transl. Ser.2, 226,
	Amer. Math. Soc., Providence, RI, 2009.
	
\bibitem{16}
	G. Wang,
	\emph{Sharp maximal inequalities for conditionally symmetric martingales and Brownian motion},
	Proc. Amer. Math. Soc. 112 (1991),
	579--586.
	
\end{thebibliography}
\end{document}